\documentclass[12pt,fleqn,emtex]{article}

\usepackage{graphicx}
\usepackage{framed}
\usepackage{times}
\usepackage{mathptmx}

\usepackage[german, american]{babel}  
\usepackage[utf8]{inputenc}

\usepackage{amsfonts}       %
\usepackage{amsmath} 
\usepackage{amssymb}        %
\usepackage{amstext}        %
\usepackage{amsthm}         
\usepackage[matrix,curve]{xypic}

\usepackage{textcomp}
\usepackage{epsfig}
\psfigdriver{dvips}
\usepackage{latexsym}
\usepackage[matrix,curve]{xypic}
\usepackage{rotating}
\rotdriver{dvips}

\newcommand{\N}{\mathbb{N}}
\newcommand{\Z}{\mathbb{Z}}
\newcommand{\Q}{\mathbb{Q}}

\newtheorem{Lemma}{Lemma}
\newtheorem{Theorem}{Theorem}
\newtheorem{Corollary}{Corollary}
\newtheorem{Proposition}{Proposition}

\DeclareMathSizes{12}{11}{8}{5} 
\usepackage{xypic,dsfont}

\begin{document}

\medskip
\centerline{\bf Vanishing theorems for constructible sheaves on}

\smallskip
\centerline{\bf abelian varieties over finite fields}

\bigskip
\centerline{(Rainer Weissauer)}

\medskip\noindent

\begin{abstract} \noindent  Let $\kappa$ be a field, finitely generated  over its prime
field, and let $k$ denote an algebraically closed field containing $\kappa$.
For a perverse $\overline\Q_\ell$-adic sheaf $K_0$ on an abelian variety $X_0$
over $\kappa$, let $K$ and $X$ denote the base field extensions of $K_0$ and $X_0$
to $k$. Then, the aim of this note is to show that  
the Euler-Poincare characteristic of the perverse sheaf $K$ on $X$ is a non-negative integer, i.e. $\chi(X,K)=\sum_\nu (-1)^\nu \dim_{\overline \Q_\ell}(H^\nu(X,K))\geq 0$. This generalizes the result of  Franecki and Kapranov
[FK] for fields of characteristic zero. Furthermore we show that $\chi(X,K)=0$ implies $K$ to be translation invariant. This result allows to considerably simplify the proof of    
the generic vanishing theorems for constructible sheaves on complex abelian varieties of [KrW].  Furthermore it extends these vanishing theorems to constructible sheaves on abelian varieties over finite fields. \end{abstract}

\bigskip
The proof in [FK] of the above estimate for the Euler-Poincare characteristic 
of perverse sheaves on abelian varieties over fields of characteristic zero as well as the results of [KrW] rely on methods from the theory of $D$-modules via the Dubson-Riemann-Roch formula for characteristic cycles. In fact, one expects a similar Riemann-Roch theorem
over fields of positive characteristic, extending the results of [AS] and generalizing the  
Grothendieck-Ogg-Shafarevich formula for the Euler-Poincare characteristic of sheaves on curves.
However, in the absence of such deep results on wild ramification we will follow a different approach
using methods of Gabber and Loeser [GL] that are based
on Ekedahl's adic formalism. 

\medskip
One of the main tools in [GL] for the study of perverse sheaves on tori
is the  Mellin transform of a perverse sheaf. This Mellin transform is a complex of sheaves on the affine spectrum of a kind of Iwasawa ring associated to the fundamental group of the torus. For us, the key property of the Mellin transform is the coherence of its cohomology sheaves, which is a consequence of the finiteness theorems of Ekedahl for adic sheaves [E]. For abelian varieties  one can define the Fourier transform of a perverse sheaf in a similar way, and most of the proofs in [GL] carry over ad hoc. Other than that, the situation for tori and abelian varieties is quite different. Whereas in the former case the non-completeness of tori leads to difficulties in defining a convolution product with nice duality properties, this is almost trivial in the case of abelian varieties. The convolution is important for the final aim to construct Tannakian categories, and both the Mellin and the Fourier transform convert the convolution product into the tensor product of complexes. 

\medskip
On the other hand, 
by fibering a torus over lower dimensional tori with one dimensional  multiplicative fibers, various important results in the torus case [GL] immediately follow from the Artin-Grothendieck vanishing theorem.
The main difficulty for abelian varieties arises from the fact that analog vanishing theorems, although needed,  
do not follow from the Artin-Grothendieck vanishing theorem in the absence of such fibrations, but have to be proven by other means. 
In [KrW]  for complex abelian varieties they are derived  from the Dubson-Riemann-Roch theorem in the case of a simple abelian variety.  In [W3], [W2] the reduction to the simple case was achieved by a tedious induction on the number of simple factors together with some characteristic $p$
argument (Cebotarev density).   In the present paper we give a new easier proof of  the following key ingredient for these vanishing theorems.

\bigskip
{\bf Main Theorem}. 
{\it For an abelian variety $X_0$ over a field $\kappa$ finitely generated field over its prime field, let  $K$ be a simple perverse sheaf on $X$ defined over $\kappa$. Then the Euler-Poincare characteristic satisfies $\chi(X,K)\geq 0$. If $\chi(X,K)=0$, then $K$ is translation invariant in the sense that $T_x^*(K)\cong K$ holds for all closed points $x$ of an abelian subvariety $A$ of $X$ of positive dimension}. 

\medskip
By the theory of monoidal perverse sheaves from [W], a degeneration argument allows to reduce the proof of this theorem to the case of finite fields $\kappa$.  
Although not stated above, then the theorem not only holds for finitely generated fields but also  for perverse sheaves on complex abelian varieties. This again can be shown by a degeneration argument, but  now via the method of Drinfeld [Dr2] involving the by now proven de Jong's conjecture. 
In the finite field case, by studying the support of the Fourier transform of monoidal perverse sheaves and by the crucial coherence properties, the proof is reduced to an elementary
statement on the existence of nontrivial one-dimensional analytic subgroups within analytic $\varphi$-subvarieties of analytic tori that are endowed with a certain automorphism $\varphi$, defined by the Frobenius. For the precise statements see proposition 4
and the related statements preceding it. That one can apply proposition 4 requires certain properties of the supports, that can be deduced from the Cebotarev density theorem [W2]. 
 
\medskip
As already shown in [KrW], the theorem stated above does imply the generic vanishing theorems for
perverse sheaves on abelian varieties, provided one also disposes over the  decomposition theorem and Hard Lefschetz theorem. To keep the range of applications flexible,
in [KrW]  certain classes ${\bf P}$ of $\overline {\mathbb Q}_\ell$-adic perverse sheaves on abelian varieties were considered that satisfy the required properties in an axiomatic way. A typical class are the perverse sheaves of geometric origin. If the abelian variety and the perverse sheaf are defined over an algebraically closed field of characteristic zero or over a finite field, from the  axiomatic treatment in [KrW] in particular one obtains the following result from the main theorem above.

\medskip
{\bf Vanishing Theorem}. {\it Let $X_0$ be an abelian variety and let $K_0$ be a $\overline {\mathbb Q}_\ell$-adic perverse sheaf on $X_0$, both defined either over the field $\kappa$ of complex numbers or a finite field. Then for most characters $\chi: \pi_1(X)\to \overline {\mathbb Q}_\ell^* $ of the etale fundamental group $\pi_1(X)$ 
the etale cohomology groups $H^\nu(X,K\otimes_{\overline {\mathbb Q}_\ell} L_\chi)$ vanish in degree $\nu\neq 0$.} 

\medskip
Here $L_\chi$ denotes the rank 1 local $\overline {\mathbb Q}_\ell$-adic system defined by the character
$\chi$. For the precise meaning of the notion \lq{most}\rq\ we refer to [KrW], but we indicate that
for absolutely simple abelian varieties $X_0$  it just means \lq{for almost all}\rq.  
      
\medskip
The vanishing theorem stated above similarly holds for perverse sheaves of geometric origin.
It might very well hold for arbitrary perverse sheaves on abelian
varieties over arbitrary algebraically closed fields. Presently this is not known, since no specialization argument as in [Dr2] for the characteristic zero case is known in the case of positive characteristic. But, with such generalizations in mind we also consider weaker vanishing statements where vanishing for most characters is replaced by the weaker statement of  vanishing for a generic character.
See theorem 2 and the discussion leading to it.

\bigskip
{\it Axiomatic setting}.
Let $k$ denote the algebraic closure of a  finite field $\kappa$ of characteristic $p$. For an abelian variety $X_0$ over $\kappa$, let $X$ be the base extension of $X_0$ from $\kappa$ to $k$
for a fixed embedding $\kappa \subset k$. Let $\Lambda$ denote ${\overline \Q}_\ell$ for some prime $\ell \neq p$. On the derived category $D_c^b(X,\Lambda)$ of $\Lambda$-adic sheaves with bounded constructible cohomology sheaves one has the convolution product $K*L = Ra_{0*}(K\boxtimes L)$ defined by the group law $a_0: X_0\times X_0\to X_0$. The convolution product makes $D_c^b(X,\Lambda)$ into a rigid symmetric monoidal category. The rigid dual $K^\vee$ of $K\in D_c^b(X,\Lambda)$ is $ (-id_X)^*D(K)$ for the Verdier dual $D(K)$ of $K$.
By definition, a character $\chi\colon \pi_1(X) \to \Lambda^*$ of the etale fundamental group $\pi_1(X)$ of $X$ is a  continuous homomorphism with values in the group of units ${\frak o}_\lambda^*$ of the ring of integers ${\frak o}_\lambda$ of a finite extension field $E_\lambda\subset \Lambda$ of $\Q_\ell$. Associated to a character $\chi$ there is a  smooth $\Lambda$-adic sheaf $L_\chi$ on $X$.
Let $\pi_1(X)_\ell$ denote the maximal pro-$\ell$ quotient of $\pi_1(X)$.
Any character $\chi$ of $\pi_1(X)$ is the product of a character $\chi_f$ of finite order prime to $\ell$, and a character
that factorizes over the pro-$\ell$ quotient $\pi_1(X)_\ell$ of $\pi_1(X)$. 

\medskip
We now fix a semisimple suspended monoidal rigid subcategory  $\bf D={\bf D}(X)$  of the triangular monoidal category $(D_c^b(X,\Lambda),*)$ that is closed under retracts and 
that satisfies the properties formulated in [KrW,$\S 5$]. In other words, we furthermore assume that:
$\bf D$ is stable under the perverse truncation functors ${}^p \tau_{\geq 0}$ and 
${}^p \tau_{\leq 0}$, and hence under the perverse cohomology functors ${}^p H^0= {}^p \tau_{\geq 0} \circ {}^p \tau_{\leq 0}$ such that $K \cong \bigoplus_{n\in \mathbb Z} {}^p H^n(K)[-n]$ holds
for ${}^p H^n(K) = {}^p H^0(K[n])$ and all $K\in {\bf D}$. The full subcategory ${\bf P}={\bf P}(X)$ of objects in ${\bf D}$ that are perverse sheaves is semisimple. 
We assume that there exist 
hard Lefschetz isomorphisms 
${}^p H^{-n}(K*L) \cong {}^p H^n(K*L)(n)$
for all $K,L$ in $\bf D$. Finally we assume that $\bf D$ is stable under twists with the character sheaves $L_\chi$. Then, for $K\in {\bf D}$ resp. $K\in {\bf P}$, the twist $K_\chi := K \otimes_\Lambda^L L_\chi$ is in ${\bf D}$ resp. in ${\bf P}$.

\medskip
An example is the category ${\bf D}$ of all $K$ in $D_c^b(X,\Lambda)$ obtained by base extension from some objects $K_0$ in $D_c^b(X_0,\Lambda)$ with the property that $K$ decomposes into a direct sum of complex shifts of irreducible perverse sheaves on $X$. Although the convolution product $*$ on ${\bf D}$, induced by the group law on $X$, makes $({\bf D},*)$ into a rigid $\Lambda$-linear monoidal symmetric category, in general the convolution product does not preserve the subcategory ${\bf P}$.  

\medskip
The case of finite fields $\kappa$ is central for this paper.
If occasionally we consider other fields $\kappa$, say that are finitely generated but not finite,
we tacitly make the same axiomatic assumptions and use analogous notation.

\medskip 
{\it Fourier transform}. We define the Fourier transform in analogy to the Mellin transform in [GL].
Since most of the arguments carry over verbatim, we restrict ourselves to give the main references
from [GL]. As in [GL, p. 509] consider the ring $\Omega_X \colon = {\frak o}_\lambda[[ \pi_1(X)_\ell ]]$, a complete noetherian local ring
 of Krull dimension $1+2\dim(X)$.
For generators $\gamma_i$ of $\pi_1(X)_\ell \cong (\Z_\ell)^{2\dim(X)}$,  this ring is 
isomorphic to the formal power series ring ${\frak o}_\lambda[[t_1,....,t_{n}]]$ in the variables 
$t_i=\gamma_i- 1$  for $n=2\dim(X)$ with coefficients in ${\frak o}_\lambda$. For ${\cal C}(X)_\ell = Spec( \Lambda \otimes_{{\frak o}_\lambda} {\frak o}_\lambda [[ \pi_1(X)_\ell]])$ as in [GL, 3.2], define the scheme ${\cal C}(X)$ as the disjoint union $\bigcup_{\chi_f} \{\chi_f\} \times {\cal C}(X)_\ell $, for $\chi_f$  running over the characters $\chi_f$ of $\pi_1(X)$ of finite order prime to $\ell$. By [GL,A.2.2.3] the closed points of ${\cal C}(X)_\ell$ are the $\Lambda$-valued points of ${\cal C}(X)_\ell$. The  $\Lambda$-valued points of the scheme ${\cal C}(X)$ can be identified with the \lq{continuous}\rq\  characters $\chi\colon \pi_1(X)\to \Lambda^*$, i.e characters in our notation.
As in loc. cit. there exists a continuous character $can_X: \pi_1(X) \to \Omega_X^*$ and an associated
local system $L_X$ on $X$, which is locally free of rank 1 over $\Omega_X$. For $K\in D_c^b(X,{\frak o}_\lambda)$
 we consider  $K\otimes_{{\frak o}_\lambda}^L L_X$ as an object in $D_c^b(X, \Omega_X)$. For the structure morphism $f\colon X \to Spec(k)$, following
 [GL, p.512 and A.1], we define the Fourier transform 
 $  {\cal F}\colon D_c^b(X,{\frak o}_\lambda)  \longrightarrow D_{coh}^b(\Omega_X) $
 by ${\cal F}(K) = Rf_*(K\otimes_{{\frak o}_\lambda}^L \Omega_X)$ analogous to the Mellin transform in loc. cit.  
By proposition A.1 of loc. cit. the functor defined by extension of scalars $- \otimes_{{\frak o}_\lambda}^L \Omega_X$
commutes with direct images for arbitrary morphisms $f:X \to Y$ between varieties $X,Y$ over $k$.
By inverting $\ell$ and passing to the direct limit over all ${\frak o}_\lambda \subset \Lambda$, we easily see that ${\cal F}$ induces a functor from ${\bf D}$ to the derived category $D_{coh}^b({\cal C}(X)_\ell)$ of ${\cal O}({\cal C}(X)_\ell)$-module sheaf complexes  with bounded coherent sheaf cohomology (see loc.cit. p. 521). The functor thus obtained
 $$ {\cal F}\colon \bigl({\bf D},*\bigr) \longrightarrow \bigl(D_{coh}^b({\cal C}(X)_\ell),\, \otimes^L_\Omega\bigr) \ $$ is a tensor functor, since ${\cal F}$ commutes with the convolution product; this follows from the arguments on p. 518 of [GL]. Similarly, $ {\cal F}\colon ({\bf D},\, *) \to (D_{coh}^b({\cal C}(X)), \otimes^L_\Omega)$ can be defined as in loc. cit.  Furthermore,
 as in [GL, cor. 3.3.2], the specialization $Li_\chi^*: D_{coh}^b({\cal C}(X)_\ell) \to D_{coh}^b(\Lambda)$, defined by the inclusion  $i_\chi\colon \{\chi\} \hookrightarrow {\cal C}(X)$ of the closed point  that corresponds to the character $\chi\in {\cal C}(X)$, has the property
 $$ Li_\chi^*\bigl({\cal F}(K)\bigr) \ =\ R\Gamma(X,K_\chi) \ .$$ 

For a complex $M$ of $R$-modules and a prime ideal ${\frak p}$ of $R$ the small support
$supp_R(M)= \{ {\frak p} \vert k({\frak p})\otimes_R^L M \not\cong 0\}$ is contained in the support
$Supp_R(M) = \{ {\frak p} \vert M_{\frak p} \not\cong 0\}$. The latter is Zariski closed in $Spec(R)$.  
For a noetherian ring $R$ and a complex $M$ of $R$-modules with bounded and coherent cohomology $R$-modules $H^\bullet(M)$
both supports coincide: $supp_R(M) = Supp_R(M)$. For the regular noetherian ring $R=\Lambda\otimes_{{\frak o}_\lambda}{\frak o}_\lambda[[ \pi_1(X)_\ell]]$ 
furthermore any object $M$ in $D_{coh}^b(R) \cong D_{coh}^b({\cal C}(X)_\ell)$ is represented by a perfect complex, i.e.
a complex of finitely generated projective $R$-modules of finite length. Notice that   
$ Li_\chi^*({\cal F}(K)) =  k({\frak p})\otimes_R^L {\cal F}(K)$ holds for the maximal ideal ${\frak p}$ of $R$
with residue field $k({\frak p})=R/{\frak p}$, defined by $\chi$.

\medskip {\it Spectrum}.
By definition, for $K\in {\bf P}$ the spectrum  ${\cal S}(K) \subseteq {\cal C}(X)(\Lambda)$ is the set of characters $\chi$ 
such that $H^\bullet(X,K_\chi)\neq H^0(X,K_\chi)$.  It is well known that the Euler characteristic does not change if $K$ is twisted by a character: $\chi(X,K_\chi)\! =\! \chi(X,K)$.  Under the assumption  $\chi(X,K)\! =\! 0$, therefore the condition $\chi \in {\cal S}(K)$ is equivalent to $H^\bullet(X,K_\chi)\neq  0$, and hence equivalent to $R\Gamma(X,K_\chi)\not\cong 0$. 

\begin{Lemma}\label{thisis1}
For $K\in {\bf P}$ with $\chi(X,K)=0$, the set  of characters ${\cal S}(K) \cap {\cal C}(X)_\ell(\Lambda)$ is the set of closed points of a Zariski closed subset
of ${\cal C}(X)_\ell$.
\end{Lemma}

{\it Proof}. For $\chi(X,K)\! =\!0$, as explained above a character $\chi\in {\cal C}(X)_\ell(\Lambda)$
is in ${\cal S}(K)$ if and only $R\Gamma(X,K_\chi) \not\cong 0$ holds, or equivalently  
if $ Li_\chi^*({\cal F}(K)) =  k({\frak p})\otimes_R^L {\cal F}(K)\neq 0$. This means $\chi\in supp_R({\cal F}(K))$. Since $R$ is noetherian and  $k({\frak p})\otimes_R^L {\cal F}(K)$ has bounded coherent cohomology, this is equivalent to the condition that $\chi$ is contained in the Zariski closed subset
$Supp_R({\cal F}(K))$. \qed

\medskip

{\it Monoidal perverse sheaves}. Monoidal objects arise as follows: Whereas the convolution of two sheaves almost never is a sheaf,   the convolution of two perverse sheaf complexes under favorable conditions is a perverse sheaf complex up to some negligible sheaf complexes [KrW2]. Therefore we mainly focus on perverse sheaves in ${\bf P}$. They define a semisimple abelian subcategory of  ${\bf D}$.
However, attached to a perverse sheaf complex $K$, the convolution product $K^\vee * K$ in general
is not a perverse sheaf complex. Hence the evaluation morphism $ eval_K:  K^\vee * K \to \delta_0$,
with values in the skyscraper sheaf $\delta_0$ at the origin, is not a morphism 
between perverse sheaves. To measure the deviation, we consider the perverse truncation
functors $^{p}\tau_{< \nu}$ on $\bf D$ to define the degree $\nu_K$ of $K$ 
to be  the largest integer $\nu$ such that $eval_K \circ  {}^p\tau_{<\nu} $ vanishes. 
Then $eval_K$ induces a nontrivial morphism of the $\nu_K$-th perverse cohomology  to $\delta_0$ $${}^p H^{\nu_K}(K^{\vee} *K)[-\nu_K] \to \delta_0\ .$$
Since $\nu_{K\oplus L}=\min(\nu_K,\nu_L)$ and $\nu_{K[n]}=\nu_K$, by our semisimplicity assumptions
one can assume that $K$ is a simple perverse sheaf and that ${}^p H^{\nu_K}(K^{\vee} *K)[-\nu_K]$ decomposes into a finite direct sum of simple perverse sheaves. Using this, it is then not difficult to show 
that for simple perverse $K$ there exists a unique simple shifted perverse direct summand $${\cal P}_K[-\nu_K]\ \subseteq \ {}^p H^{\nu_K}(K^{\vee} *K)[-\nu_K]$$ on which the morphism induced by $eval_K$ is nontrivial. See [W]. The so defined perverse sheaf ${\cal P}_K$ has particular properties, similar to those of ambidextrous objects in a rigid symmetric monoidal abelian category, and is called the monoidal perverse sheaf associated to the simple perverse sheaf $K$. 

\medskip For simple objects $K$ in ${\bf P}$ this defines  an integer in $[0,\dim(X)]$, the degree  $\nu_K$ of $K$,
and an irreducible monoidal perverse sheaf ${\cal P}_K$ in ${\bf P}$. 
By [W, lemma 1.4] the Euler-Poincare characteristic $\chi(X,K)$ of $K$ on $X$ is zero  if and only if
$\nu_K>0$; furthermore ${\cal P}_K \! \cong \! {\bf 1}$ (unit object) holds if and only if $\nu_K=0$. ${\cal P}_K$ is called a {\it monoid} in case $\nu_K>0$. We quote from [W,\! cor.\! 4, lemma\! 1] the following statement

\bigskip\noindent
{\bf Convexity lemma}.
{\it For monoids ${\cal P}_1,{\cal P}_2\in {\bf P}$ of degrees $\nu_1, \nu_2$ on an abelian variety $X$ with ${\cal P}_1 \not\cong {\cal P}_2$ the convolution product $L={\cal P}_1 * {\cal P}_2$ has degree $\nu_L > (\nu_1 +\nu_2)/2$.}

\medskip
Furthermore, one can show $\nu_K = \nu_{{\cal P}_K}$ and that ${\cal P}_K$ is the monoidal perverse sheaf also associated to ${\cal P}_K$ itself. Finally, $K$ is negligible if and only if ${\cal P}_K$ is negligible resp. if and only if $\nu_K>0$. For details see [W].

\begin{Lemma} \label{new} If $\nu_K>0$  holds for a simple object $K\in {\bf P}$,  
then ${\cal S}(K) = {\cal S}({\cal P}_K)$.
\end{Lemma}

{\it Proof}. As mentioned above, $\nu_K>0$ implies $\chi(X,K)\! =\! 0$. From the discussion of lemma 1, therefore the condition $\chi \in {\cal S}(K)$ is equivalent to $R\Gamma(X,K_\chi)=0$ 
respectively $H^\bullet(X,K_\chi)=0$. Now we use the
split monomorphisms $ K[\pm\nu_K] \hookrightarrow {\cal P}_K *K$ and 
${\cal P}_K[\pm \nu_K] \hookrightarrow K*K^\vee$ constructed in [W].
Using $(A*B)_\chi \cong A_\chi * B_\chi$, the second one implies $H^\bullet(X,({\cal P}_K)_\chi)=0$.
Hence $\nu_{\cal P_K}>0$. Indeed, we therefore see that the assertions $H^\bullet(X,K_\chi)=0$ and $H^\bullet(X,({\cal P}_K)_\chi)=0$ are equivalent. Hence ${\cal S}(K) = {\cal S}({\cal P}_K)$. \qed

\medskip 
If $\nu_{K_i}>0$ for either $i=1$ or $i=2$,  by [KrW] all simple constituents $K[n]$ of $K_1*K_2\cong \bigoplus K[n]$ satisfy $\nu_K>0$. In general, the semisimple complexes with simple constituents of vanishing Euler-Poincare characteristic define a tensor ideal ${\bf N}_{Euler}$ in ${\bf D}$. All monoids are in this tensor ideal ${\bf N}_{Euler}$. For any semisimple complex $K$ in ${\bf N}_{Euler}$, let ${\cal S}(K)$ denote the set of $\chi\in {\cal C}(X)(\Lambda)$ for which $H^\bullet(X,K_\chi)\neq 0$. Then $ {\cal S}(K \oplus K') = {\cal S}(K) \cup {\cal S}(K') $, and by the K\"unneth formula
$$ {\cal S}(K * K') = {\cal S}(K) \cap {\cal S}(K') \ $$
holds for all semisimple complexes $K, K'$ in ${\bf N}_{Euler}$.

\begin{Lemma} \label{suppo} If  for a simple perverse sheaf $K$ in ${\bf N}_{Euler} \subset {\bf D}$ and a character $\chi_f$ of order  prime to $\ell$ the scheme $Y={\cal S}(K) \cap (\{\chi_f\}\times {\cal C}(X)_\ell)$ contains an isolated closed point $y$, then this scheme has only this closed point $y$ and $K$ is a character twist
of the perverse sheaf $\delta_X:=\Lambda_X[\dim(X)]$.
\end{Lemma}

{\it Proof}. 
We may assume $\chi_f\! =\! 1$ by twisting $K$. Then $P\!:=\! {\cal F}(K)$ is a complex in $D_{coh}^b(R)$. For $\chi$ corresponding to the isolated closed point $y\in Y \subseteq {\cal C}(X)_\ell$ let $m_y$ be the associated maximal ideal of $R$ with residue field $\Lambda_y$.
Then $ R\Gamma(X,K_\chi) \cong Li_\chi^*\bigl({\cal F}(K)\bigr)  \cong P\otimes_R^L \Lambda_y$. We claim:
$H^\nu(X,K_\chi) \neq 0$ for some $\nu$ with $\vert \nu\vert \geq \dim(X)$, so $K_\chi \cong \delta_X$ follows from the known  cohomological bounds for perverse sheaves $K_\chi$. To find $\nu$, it suffices to find $i \geq j+ 2\dim(X)$ with $H^{\mu}(X,K_\chi)\neq 0$ for $\mu=i,j$.  
For that replace $R$ by its localization at $m_y$, a regular local ring of
dimension $n=2\dim(X)$.  Then $P\in D^{[a,b]}_{coh}(R)$ for some integers $a\leq b$, i.e. $H^i(P)=0$ holds for $i\notin [a,b]$  with cohomology supported in $\{ y\}$, and we may assume $H^a(P)\neq 0$ and $H^b(P)\neq 0$. For any such $Q\in D^{[a,b]}_{coh}(R)$  we claim: $H^{a-n}(Q\otimes_R^L \Lambda_y)\neq 0$ and $
H^{i}(Q\otimes_R^L \Lambda_y)= 0$ for $i\! <\! a\! -\! n$.  By decalage, using truncation triangles $\tau^{\leq a}(Q)\to Q \to \tau^{>a}(Q)\to $ and shifts, this claim is easily reduced to the case where $a=b=0$. Then $Q$ represents an $R$-module with support in $m_y$ and  $\Lambda_y$ is the only simple $R$-module with support in $m_y$. Hence by decalage the claim is reduced to the case where $Q$ is quasi-isomorphic to $\Lambda_y$. Then $Q$ is represented by the perfect minimal Koszul complex $Kos(\Lambda_y)=(0\to Q_{-n} \to ... \to Q_0 \to 0)$ with $Q_{-n} \cong Q_0 \cong R$, so our claim now follows  from $H^i(Q \otimes_R^L \Lambda_y) = H^i(Kos(\Lambda_y) \otimes_R \Lambda_y)$
with $H^{-n}(Q \otimes_R^L \Lambda_y) = Tor^R_n(\Lambda_y,\Lambda_y) \cong \Lambda_y$.
Under the assumptions on $Q$ above, also  $H^b(Q\otimes_R^L \Lambda_y)\neq 0$ and $H^i(Q\otimes_R^L \Lambda_y)=0$ 
for $i>b$. Again one reduces this to the case of the Koszul complex $Kos(\Lambda_y)$. 
Applied to $P$, this shows $H^{i}(X,K_\chi)\neq 0$ and
$H^j(X,K_\chi)\neq 0$ for $i=b$ and $j=a-2\dim(X)$.  \qed

\begin{Lemma}\label{FITW}
For an irreducible perverse sheaf $K$ on $X$, the group $\Delta_K =\{ \chi \ \vert \ K \cong K_\chi\}$
is a subgroup of the group ${\cal C}(X)(\Lambda)$ of all characters $\chi$ of $\pi_1(X)$. It is a proper subgroup unless $K$ is a skyscraper sheaf. More precisely,
let $A$ be the abelian subvariety generated by the support of the perverse sheaf $K$ in $X$ and let $K(A)$ denote the subgroup of 
characters in ${\cal C}(X)(\Lambda)$ whose restriction to $A$ becomes trivial. Then $K(A)$ is a subgroup of $\Delta_K$ and the quotient $\Delta_K/K(A)$
is a finite group.   
\end{Lemma}

{\it Proof}. Suppose $K$ is not a skyscraper sheaf. Then the support $Y$ of $K$ generates an abelian subvariety $A\neq 0$ of $X$. We may replace $X$ by this subvariety $A$. Then the natural morphism
$H^1(X,\Lambda) \to H^1(Y,\Lambda)$ is injective, and hence $\pi_1(Y,y_0) \to \pi_1(X,y_0)$
has finite cokernel [S, lemma VI.13.3, prop. VI.17.14], say of index $C$. There exists a Zariski open dense subset $U$ of $Y$
and a smooth $\Lambda$-adic sheaf $E$ on $U$, defining a $\Lambda$-adic representation $\rho$,  such that $K\vert_U \cong E[\dim(Y)]$.
Since $\rho\otimes \chi\cong \rho$ for all $\chi\in \Delta_K$, viewed as characters $\chi$ of $\pi_1(Y,y_0)$, we obtain the following bound $\#\Delta_K \leq C \cdot \dim_\Lambda(\rho)$ from the next lemma.
\qed

\begin{Lemma}\label{5}Let $\rho$ be an irreducible representation of a group $\Gamma$ on a finite-dimensional vectorspace over $\Lambda$, and let $\Delta$ be a finite group of abelian characters $\chi:\Gamma\to \Lambda^*$, defining a normal subgroup $\Gamma' = Ker(\Delta)$ such that $\Gamma/\Gamma' \cong \Delta^*$. Then $\rho\otimes \chi \cong \rho$ for all $\chi\in \Delta$ implies $\rho\cong Ind_{\Gamma'}^\Gamma(\rho')$ for some irreducible representation $\rho'$ of $\Gamma'$. In particular $$\#\Delta \ \leq\ \#\Delta \cdot \dim_\Lambda(\rho') = \dim_\Lambda(\rho)\ .$$ \end{Lemma}

{\it Proof}. For the convenience of the reader we give the proof.
If $\rho \cong Ind_{\Gamma_0}^\Gamma(\rho_0)$ for some subgroup $\Gamma' \subseteq {\Gamma_0}\subseteq
\Gamma$, we may replace the pair $(\Gamma,\rho)$ by $(\Gamma_0,\rho_0)$. Indeed, $\rho_0 \otimes (\chi\vert_{\Gamma_0}) \cong \rho_0$ for $\chi\in \Delta$ holds. To show this:
$\rho_0$ is a constituent of $Ind_{\Gamma_0}^{\Gamma}(\rho_0)\vert_{\Gamma_0} \cong \rho\vert_{\Gamma_0}$,
and therefore also a constituent of $(\rho\otimes \chi)\vert_{\Gamma_0}$. Hence $\rho_0\otimes (\chi\vert_{\Gamma_0})
\cong \rho_0^s$ by Mackey's lemma for some $s\in \Gamma$, with $s$ a priori depending on $\chi\in \Delta$. But $s\in \Gamma_0$,
since otherwise $\rho_0$ could be extended to a projective representation of $\langle \Gamma_0,s \rangle\subseteq \Gamma$, and this is easily seen to contradict the irreducibility of $\rho \cong Ind_{\Gamma_0}^\Gamma(\rho_0)$. Therefore $s\in \Gamma_0$, and this implies our claim:  $\rho_0\otimes (\chi\vert_{\Gamma_0})
\cong \rho_0$ for all $\chi\in \Delta$.

\medskip
Using induction in steps, without loss of generality we can therefore  assume   that $\rho\not\cong 
Ind_{\Gamma_0}^\Gamma(\rho_0)$ holds for any $\Gamma_0$ in $\Gamma$ such that $\Gamma' \subseteq \Gamma_0 \neq \Gamma$. We then have to show $\Gamma=\Gamma'$. If $\Gamma'\neq \Gamma$, we may now also replace the group $\Gamma'$ by some larger group $\Gamma_0$ with prime index in $\Gamma$.  Then there exists a character $\chi \in \Delta$ with kernel $\Gamma_0$. By Mackey's theorem and $\rho\not\cong 
Ind_{\Gamma_0}^\Gamma(\rho_0)$,  the restriction  $\rho\vert_{\Gamma_0}$ is an isotypic multiple $m \cdot \rho_0$ of some irreducible representation $\rho_0$ of $\Gamma_0$. Therefore $(\rho_0)^s \cong \rho_0$ holds for all $s\in \Gamma$. Hence $\rho_0$ can be extended to a representation of $\Gamma$ on the representation space of $\rho_0$ (there is no obstruction for extending the representation since $\Gamma/\Gamma_0$ is a cyclic group). By Frobenius reciprocity, this extension is then isomorphic to $\rho$; so $m=1$. In other words, the restriction of $\rho$ to $\Gamma_0$ is an irreducible representation of $\Gamma_0$, hence equal to $\rho_0$. 

\medskip
Finally, $\rho\otimes \chi \cong \rho$ implies
$\chi\hookrightarrow \rho^\vee \otimes \rho$ (as a one dimensional constituent). Therefore $\bigoplus_{\chi\in \Delta_0^*}  \chi \hookrightarrow \rho^\vee \otimes \rho$, as representations of $\Gamma$. Restricted to $\Gamma_0$, this
implies $\# \Delta_0 \cdot {\bf 1} \hookrightarrow \rho_0^\vee \otimes \rho_0$, since $\rho\vert_{\Gamma_0} \cong \rho_0$. But $Hom_{\Gamma_0}({\bf 1}, \rho_0^\vee \otimes \rho_0) \cong Hom_{\Gamma_0}(\rho_0,  \rho_0) \cong \Lambda$ since $\rho_0$ is irreducible. Hence $\# \Delta_0\! =\! [\Gamma:\Gamma_0] \! =\! 1$. This implies $\Gamma\! =\! \Gamma_0$, and hence $\Gamma\! =\! \Gamma'$. \qed

\begin{Proposition}\label{smallsupport} Suppose $\dim(X)>0$.
Then for any finite set $\{{\cal P}_1,...,{\cal P}_m\}$ of monoids  in $\bf P$, there exist characters $\chi\in {\cal C}(X)_\ell$  such that $\chi \not\in \bigcup_{i=1}^m {\cal S}({\cal P}_i)$.
\end{Proposition}

{\it Proof}. 
Since  the spectrum of $R=\Lambda\otimes_{{\frak o}_{\lambda}} {{\frak o}_{\lambda}}[[x_1,...,x_n]]$ is not the union of finitely many Zariski closed proper subsets for $n\! =\! 2\dim(X)>0$,  it suffices that the spectrum ${\cal S}({\cal P})_\ell\! =\! {\cal S}({\cal P}) \cap {\cal C}(X)_\ell(\Lambda)$ of each monoid ${\cal P}$ is the set of closed points of some proper Zariski closed subset of ${\cal C}(X)_\ell$.  
We prove this by descending induction on the degree $\nu_{\cal P}$.
For $\nu_{\cal P}\! =\! \dim(X)$ this is clear, since in this case  ${\cal S}({\cal P})$ is a single point ([W,\! lemma\! 1]).
For a given monoid  ${\cal P}$ and fixed $\nu\! =\! \nu_{\cal P}\! < \! \dim(X)$, assume our assertion is true for all monoids ${\cal Q}$ of degree $ \nu_{\cal Q} > \nu$. 
By lemma \ref{FITW} there exists a character $\chi\in {\cal C}(X)_\ell$ such that ${\cal P}_\chi\not\cong {\cal P}$. Since ${\cal P}$ and ${\cal P}_\chi$ have the same degree $\nu=\nu_{\cal P}$, this implies that all constituents $K[m], K\in {\bf P}$ of ${\cal P} * {\cal P}_\chi$ have associated monoids ${\cal P}_K$ of degree $> (\nu_{\cal P} + \nu_{{\cal P}_\chi})/2 = \nu$ by the convexity lemma. Hence ${\cal S}_\ell({\cal P} * {\cal P}_\chi)$ is contained in a proper Zariski closed subset of the spectrum ${\cal C}(X)_\ell$, by  lemma \ref{new} and the induction assumption.
Suppose ${\cal S}({\cal P})_\ell$ were not contained in a proper Zariski closed subset of ${\cal C}(X)_\ell$.
Then ${\cal S}({\cal P})_\ell = {\cal C}(X)_\ell(\Lambda)$, and therefore ${\cal S}_\ell({\cal P}_\chi) = {\cal S}_\ell({\cal P}) \cap {\cal S}_\ell({\cal P}_\chi)$. Hence ${\cal S}({\cal P}_\chi)_\ell$ would be contained in a proper Zariski closed subset of ${\cal C}(X)_\ell$. Indeed, this would  follow from ${\cal S}_\ell({\cal P}_\chi) = {\cal S}_\ell({\cal P}) \cap {\cal S}_\ell({\cal P}_\chi) =  {\cal S}_\ell({\cal P} * {\cal P}_\chi)$ and the induction assumption. 
On the other hand, ${\cal S}({\cal P}_\chi)_\ell = \chi^{-1} \cdot {\cal S}({\cal P})_\ell = {\cal C}(X)_\ell(\Lambda)$. This gives a contradiction, and proves our claim for the fixed  degree $\nu$. Now proceed by induction.
\qed

\medskip
For $K\in {\bf P}$ the $\ell$-spectra ${\cal S}(K_{\chi_f^{-1}})_{\ell} = {\cal S}({K}) \cap (\{\chi_f\}\times {\cal C}(X)_\ell)(\Lambda) \subseteq {\cal S}(K)$ at some  given point $\chi_f$ of ${\cal S}(K)$ are the $\Lambda$-valued points of a Zariski closed subset
of $\{\chi_f \}\times {\cal C}(X)_\ell$ by lemma \ref{thisis1}. Replacing $K$ by $K_{\chi_f}$ we may always assume $\chi_f=1$.

\begin{Corollary} \label{existence} For
any semisimple complex $K\in {\bf D}$ contained in ${\bf N}_{Euler}$, there exists in ${\cal C}(X)_\ell(\Lambda)$ 
a character $\chi \notin {\cal S}(K)$.
\end{Corollary}

{\it Proof}.  Since ${\cal S}(K)={\cal S}({\cal P}_K)$
for simple $K$ and ${\cal S}(\bigoplus_{i=1}^m K_i[n_i]) \subseteq \bigcup_{i=1}^m {\cal S}(K_i)$, this is an immediate consequence of lemma \ref{new} and proposition \ref{smallsupport}. \qed

\begin{Theorem}\label{thm1}
For arbitrary $K\in {\bf P}$, the Euler-Poincare characteristic $\chi(X,K)$ is non-negative. Hence, in particular, the reductive supergroup ${\bf G}(K)$ attached to $K$ in [KrW,$\S 7$] is a reductive algebraic group over $\Lambda$. 
\end{Theorem}

{\it Proof}. We may assume that $K$ is irreducible.  Then, to show $\chi(X,K)\geq 0$,
it is enough to show the existence of a character $\chi$ such that $H^\nu(X, K_\chi)=0$
holds for all $\nu\neq 0$. Then  $\chi(X,K)=\chi(X,K_\chi)= \dim_\Lambda(H^0(X,K_\chi))$, and the claim obviously follows 
from $\dim_\Lambda(H^0(X,K_\chi))\geq 0$. So, we have to find a character $\chi\notin {\cal S}(K)$.
By [KrW,$\S 9$], for all irreducible perverse sheaves $K$ there exists a perverse sheaf $T$ in ${\bf N}_{Euler}$, depending on $K$,  such that $H^\bullet(X, K_\chi)\neq H^0(X,K_\chi)$ holds if and only if $\chi \in {\cal S}(T)$. Hence, by corollary \ref{existence} there exists a character $\chi\notin {\cal S}(T) = {\cal S}(K)$. 
\qed

\medskip The crucial fact that ${\cal S}(K)$ is the spectrum ${\cal S}(T)$ for an object $T$ in ${\bf N}_{Euler}$, already exploited in the proof of the last theorem,
furthermore
implies 

\begin{Theorem} \label{subset}
For any $K\in {\bf P}$ on  $X$ and any character $\chi_f$ of $\pi_1(X)$ of order prime to $\ell$, the
set of characters $\chi\in {\cal C}(X)_\ell(\Lambda)$ for which $\chi_f\chi$ is in $ {\cal S}(K)$ is the set  of closed points of a proper Zariski closed subset of ${\cal C}(X)_\ell$.  
\end{Theorem}

\medskip
For base fields $F$ of characteristic $p>0$, the following corollary now easily  follows from theorem \ref{thm1} by a specialization argument.
For the case of fields $F$ of characteristic zero see [FK]; but our argument could also be extended to the characteristic zero case.

\begin{Corollary} \label{c2}
For $\overline\Q_\ell$-adic  perverse sheaves $K_0$ on abelian varieties $X_0$ defined over a field $F$ finitely generated over its prime field, with base extensions $K$ resp. $X$ to an algebraic closure of $F$, 
the Euler-Poincare characteristic $\chi(X,K)$ is non-negative.  
\end{Corollary}

\medskip
{\it Translation invariance}. An irreducible perverse sheaf $K$ on an abelian
abelian variety $X$ will be called {\it translation invariant} if there exists an abelian subvariety $A$
of $X$ of dimension $\dim(A)>0$ such that one of the two equivalent conditions of the next lemma holds.

\begin{Lemma} \label{equiv} For an irreducible perverse sheaf $K$ on $X$ and an abelian subvariety $A$ of $X$
the following assertions are equivalent: a) $T_x^*(K)\cong K$ holds for all closed points of $X$ of $A$, resp. b) $$K\cong  L_\chi \otimes_{\Lambda_X} q^*(\tilde K)[dim(A)]$$ holds for the quotient morphism $q : X \to \tilde X = X/A$, some 
perverse sheaf $\tilde K$ on $\tilde X$ and some character $\chi: \pi_1(X, 0) \to \Lambda^*$.
\end{Lemma}

\medskip 
{\it Proof}. $K$ corresponds to a $\Lambda$-adic representation
$\phi: \pi_1(U) \to Gl(n,\Lambda)$, for a suitable open Zariski dense subset $U$ of the support $Z$ of $K$. If $K$ is invariant under translations by the closed points in $A$, the open subset $U$ can be  chosen
such that $U+A=U$. For $\tilde U = U/A$ the projective morphism $\pi_1(q\vert_U)$ induces an exact sequence $\pi_1(A) \to \pi_1(U) \to \pi_1(\tilde U)\to 0$. Its first morphism $\sigma$ is injective, since the composition $\rho\circ \sigma$ with $\rho\!:\! \pi_1(U) \to \pi_1(X)$ induced from  $U\!\hookrightarrow \! X$ is the injective morphism $\pi_1(A) \to \pi_1(X)$ that is obtained from the inclusion $A\hookrightarrow X$. 
Hence $\pi_1(A)$ is a normal subgroup of $\pi_1(U)$.  We claim that $\pi_1(A)$ is in the center of $\pi_1(U)$. 
Indeed, for $\alpha\in \pi_1(A)$ and $\gamma\in \pi_1(U)$ there exists an $\alpha'\in \pi_1(A)$ such that
 $\gamma \alpha \gamma^{-1} = \alpha'$, and hence $\rho(\gamma) \rho(\alpha) \rho(\gamma)^{-1}\! =\! \rho(\alpha')$.  Since
$\pi_1(X)\! =\! H_1(X)$ is abelian, therefore $\rho(\alpha)\! =\! \rho(\alpha')$. Because $\rho\circ \sigma$ is injective, hence $\alpha=\alpha'$.  Now, since $\pi_1(A)$ is a central subgroup of $\pi_1(U)$,  there exists a character $\chi$ of $\pi_1(A)$ such that $\phi(\alpha\gamma)=\chi(\alpha)\phi(\gamma)$ holds for the irreducible representation $\phi$ of $\pi_1(U)$, for $\alpha\in \pi_1(A)$ and $\gamma\in \pi_1(U)$. Since $\pi_1(\tilde X)$ is a free $\Z_\ell$-module, any character $\chi$ of $\pi_1(A)$ with values in $\Lambda^*$ can be extended to a character $\chi_X$ of $\pi_1(X)$. Thus $\chi_X^{-1} \otimes\phi$ is an irreducible representation, which is trivial on $\pi_1(U)$; so it defines an irreducible representation of $\tilde U$ and an irreducible perverse sheaf $\tilde K$ on $\tilde U$
such that $L_{\chi_X}^{-1}\otimes K = q^*(\tilde K)[\dim(A)]$ holds on $U$. 
Let $\tilde K$ also denote the intermediate extension of $\tilde K$ to $\tilde Z$, which is an irreducible perverse sheaf on $\tilde Z$. Since $q:Z\to \tilde Z$ is a smooth morphism with connected fibers, $q^*[dim(A)]$ is a fully faithful functor from the category of perverse sheaves on $\tilde Z$ to the category of perverse sheaves on $Z$. Hence $L=q^*(\tilde K)[\dim(A)]$ is an irreducible perverse sheaf on $Z$. Since $K$ and $L$ are irreducible perverse sheaves on $Z$ whose restrictions on $U$ are isomorphic, 
the perverse sheaves $K$ and $L$ are isomorphic. This proves the nontrivial direction. \qed

\medskip
 The characterization b) of lemma \ref{equiv} implies $\chi(X,K)=0$
 for translation invariant irreducible perverse sheaves $K$.
Over finite fields we have the following converse.

\begin{Theorem} \label{3}
For an abelian variety $X_0$ over a finite field $\kappa$, let  $K$ be a simple perverse sheaf on $X$ defined over $\kappa$. If  $\chi(X,K)=0$, then $K$ is translation invariant. 
\end{Theorem}

\medskip
Before we start with the proof of theorem 3, let us gather a couple of remarks.
These remarks reduce the assertion of theorem 3 to elementary statements 
on the geometry of nonarchimedian analytic varieties endowed with a Frobenius action, as formulated in proposition 3 and 4. The key steps for this reduction to analytic geometry are explained in remark 2 and 3 below centering around proposition 2. In fact, theorem 3 is thus reduced to proposition 2, which in turn will follow from proposition 3 resp. the related proposition 4. 

\medskip
{\it Preliminary remarks}. 1)  
If  the perverse sheaf $K$ is not translation invariant under the conditions of theorem \ref{3}, ${\cal S}(K)$ must contain infinitely many torsion characters by [W2, lemma 16]. These torsion characters are defined over finite extension fields of $\kappa$. So we can replace $K$ by a suitable twist with some torsion character $\chi$, thereby passing to a suitable finite extension field of $\kappa$ if necessary, so that
${\cal S}_\ell(K) ={\cal S}(K) \cap {\cal C}(X)_\ell(\Lambda)$ is nonempty and contains the trivial character $\chi_0$.
The character $\chi_0$ is a fixed point of Frobenius $F_\kappa$, and by lemma \ref{suppo} we may suppose that $\chi_0$ is not an isolated point of the scheme ${\cal S}_\ell(K)$. 
 
\medskip 
Remark 2) We may assume that $K$ is a monoidal perverse sheaf. Recall the monoidal perverse sheaf ${\cal P}_K$ associated to an irreducible perverse sheaf $K$ on $X$ and the corresponding {\it degree} $\nu_K$. By definition, up to a degree shift this perverse sheaf is a direct summand of $K^\vee*K$ for the convolution product $*$ of $K$ with its Tannaka dual $K^\vee = (-id_X)^*(DK)$. Since ${\cal P}_K$ appears with multiplicity one in ${}^p H^{-\nu_K}(K^\vee*K)$ by [W, lemma 1.6], it defines a pure Weil perverse sheaf  on $X$ coming from a perverse sheaf ${\cal P}_{K_0}$ over $\kappa$. Hence we can assume $K_0(\alpha)={\cal P}_{K_0}$ for a generalized Tate twist $K_0(\alpha)$ of $K_0$, defined by some $\alpha=\alpha(K_0)$ in $\Lambda^*$.
Recall, $\chi(X,K)=0$ is equivalent to $0<\nu_K$. Since $\nu_K = \nu_{\cal P_K}$ holds by [W, lemma 3], this implies
$\chi(K,{\cal P}_K)=0$.  Furthermore ${\cal P}_K$ is translation invariant if and only if $K$ is translation invariant; see [W, lemma 2.3]. So from now on suppose that $K$ is a monoidal perverse sheaf on $X$. To prove theorem \ref{3} it then suffices to
show

\begin{Proposition} \label{mono}
For an abelian variety $X_0$ over a finite field $\kappa$, any  irreducible monoidal perverse sheaf $K$ on $X$ defined over $\kappa$ is of the form $$  K \ \cong \  \delta_A \otimes_{\Lambda_X} L_\chi \ ,$$
for an abelian subvariety $A$ of $X$ so that $\nu_K =\dim(A)$ and a character $\chi: \pi_1(X)\to \Lambda^*$ of finite order. Here $\delta_A$ denotes the constant perverse sheaf  with support on $A$.
\end{Proposition}

For the proof of proposition \ref{mono}, we may normalize $K$ as in remark 1).
We may furthermore replace $X$ by the abelian subvariety $A$ of $X$, generated by the support of the monoidal perverse sheaf $K$ in $X$. Then, by lemma \ref{FITW} we have
$$ \Delta_K \colon =\{ \chi \ \vert \ K \cong K_\chi\} < \infty  \ .$$
Now assuming $X=A$, it suffices  to show that $\Delta_K < \infty$ implies $\nu_K=\dim(X)$. 
Recall, by [W1, lemma 1.2], for an irreducible monoidal perverse sheaf $\nu_K=\dim(X)$ resp. $\nu_K=0$ holds if and only if $K \cong L_\chi$ holds for some character $\chi$ of the fundamental group (i.e. translation invariance by $X$) resp. $K\cong \delta_0$ holds. In particular, $\nu_K=0$ corresponds to the trivial case where the dimension $\dim(A)$ of the support is zero.  
So, if proposition \ref{mono} does not hold, there would exist a counterexample with finite $\Delta_K$ for which $0<\nu_K < \dim(X)$ holds.  So, for the proof we will suppose that $K$ is such a critical monoidal counterexample
for which $\nu_K$ is chosen maximal with respect to the property $0<\nu_K < \dim(X)$, and then argue by contradiction.

\medskip
Remark 3) Since $\Delta_K < \infty$, we may assume that
$$ \mbox{ \it ${\cal S}_\ell(K)$ does not contain an infinite abstract 
subgroup $A$ of characters in ${\cal C}(X)_\ell(\Lambda)$.} \ $$ 
Indeed, otherwise we immediately get a contradiction that proves proposition \ref{mono}. This is shown as follows: Since  
$A$ is a group, for $\chi\in A$ we get $\chi^{-1}A \cap A = A$ and hence
$$   A \subseteq  {\cal S}(K_{\chi})  \cap  {\cal S}(K) \ .$$
Since $K_{\chi} \not\cong K$
for almost all $\chi\in A \subset {\cal S}(K)_\ell $, on the other hand
 by  the maximality of the degree $\nu_K$ $$ {\cal S}(K_{\chi})  \cap  {\cal S}(K) =  {\cal S}(K_{\chi} *K)$$ is a finite set for almost all $\chi\in A$ by the convexity lemma. A contradiction. 

\goodbreak 
\medskip
Remark 4) Using a theorem of Drinfeld [Dr] and its variant for perverse sheaves [W2, thm. 5], for our subsequent proof of proposition \ref{mono}
we may replace the $\Lambda$-adic perverse sheaf $K$ by a $\Lambda'$-adic perverse sheaf $K'$ for the algebraic closure $\Lambda'$
of $\Q_{\ell'}$ for some other prime $\ell'$ different from $p$, without changing the Euler characteristic: $\chi(X,K')=\chi(X,K)=0$. In fact, by the Grothendieck-Lefschetz etale fixed point formula the vanishing of the Euler characteristic of $K$ can be expressed  in terms the functions $f_m^K(x)$ for $x\in X(\kappa_m)$ attached to the perverse sheaf $K_0$ on $X_0$ in [W2], for all finite extension fields $\kappa_m$ of $\kappa$.  Similarly, translation invariance of $K$ resp. $K'$ can also be expressed in terms of the functions $f^K_m(x)$ resp. $f_m^{K'}(x)$. This justifies our claim, since $\tau(f^K_m(x))=f_m^{K'}(x)$ holds for a suitable field isomorphism $\tau: \Lambda \cong \Lambda'$.

\medskip
Remark 5)  The Frobenius substitution $F_\kappa$ acts
continuously on the Tate module $H_1(X,\Z_\ell)$ of the abelian variety $X$.
Choosing a basis, we identify $H_1(X,\Z_\ell)$
with $\Z_\ell^n$.  The associated $\Q_\ell$-linear map $F_\kappa$ on $H_1(X,\Q_\ell)$ can be diagonalized
over some finite dimensional extension field $E_\lambda$ of $\Q_\ell$.
Its eigenvalues $\alpha_1,...,\alpha_n\in E_\lambda^*$ are Weil numbers of weight $w>0$, i.e. they are algebraic numbers and units at all nonarchimedean places outside $p=char(k)$ with absolute value $p^{w/2}$ at the archimedean places and $E=\Q(\alpha_1,...,\alpha_n)$ is a finite extension of $\Q$.
So we can choose rational primes $\ell'$ (for some large prime $\ell' \neq p$) in such a way  that
this extension splits over $\ell'$,  thereby changing the coefficient system $\Lambda$ of $K$ into $\Lambda'$ via remark 4 from above. Therefore, without restriction of generality, we can assume $\alpha_1,...,\alpha_n \in \Z_\ell$  so that also
$F_\kappa\in Gl(n,\Z_\ell)$ is diagonalized by conjugation within the group $Gl(n,\Z_\ell)$; look at the diagonalizing matrix in $Gl(n,E)$ and
choose $\ell'$ large enough. 

\medskip
Remark 6)  We identified ${\cal C}(X)_\ell$ with the spectrum of the ring $R_n = \Lambda \otimes_{\Z_\ell} {\Z}_\ell[[t_1,....,t_{n}]]$,
viewing its closed points as characters of $\pi_1(X)$ with values in $\Lambda^*$. Thus we identify the closed points ${\cal S}_\ell(K)(\Lambda)$ of ${\cal S}_\ell(K)$
with characters in ${\cal C}(X)_\ell(\Lambda)$. By remark 5)  we can choose $\ell$ and
coordinates $t_1,...,t_n$ such that $F_\kappa (t_i) = \alpha_i \cdot t_i + P_i(t_i)$ holds for power series
$P_i$ with leading degree $\geq 2$. Indeed $$F_\kappa(\log(1+t_i)) = \alpha_i \cdot \log(1+t_i) \ .$$
Notice, these are not power series in $R_n$.  Instead we have to pass to the subring $A_n$ of locally convergent power series (see the appendix) in the formal completion $\hat R_n = \Lambda[[t_1,....,t_n]]$ of the ring 
$R_n$ with respect to the maximal ideal $(t_1,...,t_n)\subset R_n$. In $A_n$ we dispose
over new local parameters $x_i = \log(1+ t_i)$ of $A_n$. Since $\gamma_i = 1 +t_i$ are topological generators of $\pi_1(X)$, the natural action
of $Gl(n,\Z_\ell)$ on $R_n$ carries over to a continuous action of $Gl(n,\Z_\ell)$ on $A_n \subset \hat R_n$ and $\hat R_n$, that is linear on the local parameters  $x_1,...,x_n$ of $A_n$  given by  $x_i \mapsto \sum_{ij} \varphi_{ij}x_j$, provided $\varphi$ is given by the matrix $(\varphi_{ij}) \in Gl(n,\Z_\ell)$. By remark 5 we can therefore assume that the endomorphism $F_\kappa$ acts diagonally so that for $i=1,...,n$ the induced continuous ring automorphism $\varphi$ of $A_n$ is given by $\varphi(x_i)=\alpha_i \cdot x_i $ on the local parameters $x_i$. 

\medskip
{\it Proof of the main theorem}. The last theorem 3 together with corollary 2 imply the main theorem as stated in the introduction in the case of finite fields. Before we proceed with the proof of theorem 3, let us briefly indicate how
one extends the main theorem from the finite field case to the case of finitely generated fields. 
Concerning the assertion of theorem 3, as above it suffices to prove the analog of proposition 2 over finitely generated fields $\kappa$. If it were not true,  we claim there exists a monoidal perverse sheaf
${\cal P}$ defined over $\kappa$ such that ${\cal H}^{-\nu_{\cal P}}({\cal P})$ is a skyscraper sheaf
and such that $\nu_{\cal P} \neq 0$. Indeed, such a monoidal perverse counterexample descends to a monoidal perverse counterexample $\tilde{\cal P}$ on $\tilde X= X/Stab({\cal P})$. By this, we may assume that $Stab({\cal P})$ is trivial, and hence that 
${\cal H}^{-\nu_{\cal P}}({\cal P}) \cong \delta_0$
holds,  using [W, lemma 1.4] and lemma 6. We need to show that this implies ${\cal P}\cong \delta_0$. For this choose a family of abelian varieties over a basis of finite type whose generic fiber is $X_0$. If we replace the basis by a suitable Zariski open dense subset $U$, ${\cal P}$ extends to a relative perverse sheaf  over $U$. If $U$ is suitably chosen,
by generic base change ([BBD,p.154 -157], and [KrW2, appendix] in the present context) we can specialize to a closed point in general position with a finite residue field. We can achieve that the reduction $\overline {\cal P}$ remains a perverse monoidal sheaf on the reduction $\overline X$ of $X$ and satisfies $\nu_{\overline{\cal P}} = \nu_{\cal P}$ and ${\cal H}^{-\nu_{\overline{\cal P}}}(\overline{\cal P}) \cong \delta_0$. Indeed, if $U$ is suitably chosen, the spezialization base change commutes with the convolution product and Tannaka duality. 
Now, ${\cal H}^{-\nu_{\overline{\cal P}}}(\overline{\cal P}) \cong \delta_0$ implies that the stabilizer of $\overline{\cal P}$ has dimension zero. Furthermore, the specialization $eval_{\overline{\cal P}}=0$ of the evaluation morphism $eval_{{\cal P}}=0$ is not zero since otherwise ${\overline{\cal P}}$ would be zero. This implies $\nu_{\overline{\cal P}}=
\nu_{\cal P}$. By [W, lemma 1.4] this information suffices to see that ${\overline{\cal P}}$ is a monoidal perverse sheaf.  Since the stabilizer of $\overline{\cal P}$ is trivial,  theorem 3 therefore implies $\nu_{\overline{\cal P}}=0$. Hence $\nu_{\cal P}=0$, and therefore ${\cal P}\cong \delta_0$.
\qed

\medskip
This being said and with the above preliminary remarks in mind, let us start to explain the strategy for the proof of theorem 3.
As already explained in remark 2, for theorem 3 it suffices to prove proposition 2.
 
\medskip
{\it Strategy of proof}. We prove proposition \ref{mono} by contradiction. 
In proposition \ref{invariant} below we construct  an infinite subgroup $A$ of the spectrum ${\cal S}_\ell(K)$
and thus derive the desired contradiction from remark 2 and remark 3 above.
As a  
subscheme of ${\cal C}(X)_\ell$, the spectrum ${\cal S}_\ell(K)$ of $K$ defines an ideal of the ring
$R_n\cong \Lambda\otimes_{{\frak o}_\lambda} {\frak o}_\lambda[[ \pi_1(X)_\ell]]$.
The Frobenius automorphism $F_\kappa$ 
acts  on  $\pi_1(X)_\ell \cong H_1(X,\Z_\ell)$. This action 
extends to a continuous  action of $F_\kappa$ by ring automorphisms of $R_n$. 
Since the perverse sheaf $K$ is defined over $\kappa$, 
its support ${\cal S}_\ell(K)$ is invariant under the Frobenius automorphism on ${\cal C}(X)_\ell$. 
Enlarging the finite ground field $\kappa$, hence replacing the Frobenius by a suitable power, we may assume
that the Frobenius stabilizes each irreducible component of ${\cal S}_\ell(K)$, so in particular stabilizes
the irreducible component of ${\cal S}_\ell(K)$ that contains the trivial character. This irreducible component $Spec(R_n/I)$, defined by an ideal $I \subset R_n$, has Krull dimension 
$\dim(R_n/I) \geq 1$ by remark 1 above.
Since the trivial character corresponds to the maximal ideal $m=(t_1,...,t_n)$ of $R_n$, we get
$I \subsetneqq m$. By theorem \ref{subset} we also know $I\neq 0$. Therefore $$  
\{0\} \subsetneqq  I \subsetneqq m \ .$$
We say that an ideal $J\subset R_n$ is defined over $\kappa$ 
if $F_\kappa(J) = J$ holds, i.e. if it is stable under the action of the Frobenius $F_\kappa$.
Replacing the component by its reduced subscheme, we may suppose that $I$ is a prime ideal of $R_n$
defined over $\kappa$.  Then we apply the next proposition to complete the proof of proposition \ref{mono} and theorem \ref{3}.

\begin{Proposition} \label{invariant} {For any prime ideal $I\subset R_n$ defined over $\kappa$ that contains the trivial character in its spectrum, either $I$ is the maximal ideal $m$, or 
there exists an infinite subset $A$ in the $\Lambda$-valued closed points of its spectrum, so that $A$ defines an abstract group when viewed as a subset of the group of characters ${\cal C}(X)_\ell(\Lambda)$.}
\end{Proposition}

\medskip
We start with some general remarks. We assumed $I\subseteq m$. The ideals $I \subseteq m \subset R_n$ correspond one-to-one to ideals in the localization $\tilde R_n=(R_n)_{m}$ of $R_n$ at the maximal ideal $m=(t_1,..,t_n)$. The completion
of $R_n$ or $\tilde R_n$ at the maximal ideal $m$ is the ring of formal power series in the variables $t_1,...,t_n$
over the field $\Lambda$ that contains the subring $A_n=\Lambda\{\!\{ t_1,..,t_n\}\!\}$ of locally convergent power series in $t_1,...,t_n$, so that $$R_n \subset \tilde R_n \subset A_n \subset \hat R_n\ .$$ Since the power series $\log$ and $\exp$ are locally convergent,  the coordinate substutions $x_i = \log(1+t_i)$ and $1+t_i = \exp(x_i)$ for $i=1,...,n$ define an automorphism of $A_n$. Hence
$A_n = \Lambda\{\!\{x_1,...,x_n\}\!\}$. 
 The ring automorphism $F_\kappa: R_n \to R_n$, 
corresponding to the Weil numbers $\alpha_1,...,\alpha_n \in \Lambda^*$ of weight $w>0$,
extends to the completion $\hat R_n$ of $R_n$ and preserves the subrings $\tilde R_n$ and $A_n$.
In the coordinates $x_i$ it induces the continuous ring automorphism $\varphi: A_n \to A_n$ (and also of $\hat R_n$) defined by $$\varphi(x_i)=\alpha_i \cdot x_i \ .$$  We say an ideal $J$ of $A_n$ is a $\varphi$-ideal if $\varphi(J)=J$ holds.  An ideal $J \subseteq m \subset R_n$ defined over $\kappa$
extends by completion to a $\varphi$-ideal $\hat J$ of $\hat R_n$, or $\hat J \cap A_n$ of $A_n$.
For what follows, notice these easy facts: The radical ideal of a $\varphi$-ideal is a $\varphi$-ideal.
For a $\varphi$-ideal $J$ in $A_n$ its prime components are permuted by $\varphi$.
If we replace the automorphism $\varphi$ by a suitable power $\varphi^k$, without restriction of generality we may assume  that $\varphi$
stabilizes all irreducible components of $J$.

\medskip
{\it Proof of proposition \ref{invariant}}. Suppose $I$ is generated by power series $P_j(T), j\in M$ as an ideal in $R_n$. Then, as explained in the appendix, the $\varphi$-ideals $\hat I$ in $\hat R_n$ resp. $\hat I \cap A_n$ in $A_n$ are generated  by the locally convergent power series $p_j(X):=p_j(x_1,...,x_n)
= P_j(\exp(x_1)-1,..., \exp(x_n)-1)$ for $j\in M$, and $I$ is maximal in $R_n$ iff $\hat I \cap A_n$ is maximal in $A_n$ iff $\hat I$ is maximal in $\hat R_n$.
For the proof we can assume $I\neq m$  without restriction  of generality, hence $\dim(\hat A_n/(\hat I \cap A_n)>0 $ holds for the $\varphi$-ideal $\hat I \cap A_n$. $R_n$ is noetherian, hence we can assume that $M$ is finite.
By the next proposition, there exist $\lambda_1,...,\lambda_n \in \Lambda$ so that all $p_j(\lambda_1\cdot x,...,\lambda_n \cdot x)\in \Lambda$ are  convergent and vanish for $j\in M$ and all $x\in \Lambda$ 
sufficiently small, i.e. $\vert x \vert < \varepsilon $ for some $\varepsilon$. Indeed, since also the ring $A_n$ is noetherian,
the $\varphi$-ideal $\hat I \cap A_n$ has only finitely many prime components. By replacing  $\hat I \cap A_n$
with one of these prime components and $\varphi$ by a finite power of $\varphi$, the assumptions of proposition \ref{P} 
are satisfied. For $x\in \Lambda$ of sufficiently small absolute value $\vert x\vert < \varepsilon$ the exponential map
$(\exp(\lambda_1 \cdot x),..., \exp(\lambda_n\cdot x))$ converges and has values  contained in the open poly disk $\{ (t_1,...,t_n) \in \Lambda^n \  \vert \ \vert t_i\vert < 1 \mbox{ for all}  \ i=1,...,n \}$ whose points are the $\Lambda$-values closed points of $R_n$, and allows to defines the desired infinite abstract subgroup $A$ of points  with the property $A \subset Spec(R_n/I)(\Lambda) \subseteq Spec(R_n)(\Lambda)$. \qed

\begin{Proposition}\label{P}
For a prime $\varphi$-ideal $I$ contained in the maximal ideal $m$ of $A_n$, but different from $m$,
there exist $\lambda_1,...,\lambda_n\in \Lambda$ not all zero, so that $I$ is in the kernel of 
the evaluation homomorphism $\pi: A_n\twoheadrightarrow A_1=\Lambda\{\!\{ x\}\!\}$ defined by $x_i \mapsto \lambda_i \cdot x$
for $i=1,...,n$ such that the following compatibility conditions hold: After replacing $\varphi$ by a suitable power, $\alpha:=\alpha_i $ does not depend on $i$ for all indices $i$ for which $\lambda_i\neq 0$ holds. 
\end{Proposition}

{\it Proof}. We prove this for the one-dimensional cases $\dim(A_n/I)=1$ first. We remark that if $\pi$ exists in a one-dimensional case, then $I=Kern(\pi)$ by dimension reasons.
Since the case $n=1$ is trivial, assume $n=2$.  Then $I\neq\{0\}$ and $I\neq m$, and  $0\subsetneqq I$ implies $ht(I)\geq 1$ and $I \subsetneqq m$
implies $\dim(A_2/I) \geq 1$. Because for prime ideals $I$ the inequality $ht(I) + \dim(A_2/I) \leq \dim(A_2)$ holds, hence  $I$ is a prime ideal of height 1. Since $A_2$ is factorial (lemma \ref{WS}), the prime ideal $I$ is a principal ideal of $A_2$  (lemma \ref{comm}) generated by some power series $f(x_1,x_2)\in A_2$. Since $I$ is a $\varphi$-ideal,
$$ f(\alpha_1\cdot x_1, \alpha_2 \cdot x) = u(x_1,x_2) \cdot f(x_1, x_2)$$ holds for some unit $u(x_1,x_2) \in A_2$.
In the completion $\hat R_2$ of $A_2$ we can find another unit $h(X)\in 1 + (x_1,...,x_n) \subset \hat R_2$ and some constant $c\in \Lambda^*$ so that
$u = c\cdot h/h^\varphi  $ holds.  We construct this unit $h(X)$ in $\hat R_2$ for $h_\nu = \prod_{i=1}^\nu(1+z_i)$ with $z_i \in m^i$ as
an infinite  product $h = \lim_{\nu\to \infty} h_\nu$,
requiring $h_\nu/h_\nu^\varphi = uc^{-1} + y_\nu$ for $ y_\nu \in m^{\nu +1}$. We define $z_\nu$ recursively by
solving the equation $\varphi(z_\nu) - z_\nu = c y_{\nu -1} u^{-1} $ modulo $m^{\nu+1}$, using that
$\varphi - id$ acts on $m^\nu/m^{\nu +1}$ by an invertible $\Lambda$-linear map. So the formal power series $g(X)=f(X) \cdot h(X)$ is invariant under $\varphi$, in the sense that the following identity
$g(\alpha_1 \cdot x_1,\alpha_2 \cdot x_2) = c \cdot g(x_1,x_2)$ holds. This
implies that $c$ is finite product of $\alpha_1$ and $\alpha_2$
and hence a Weil number of weight $k \cdot w$ for some integer $k$, so $g(x_1,x_2)$ must be a homogenous polynomial of degree $k$ in $x_1,x_2$. The ideals $(f(X)) \subset A_2$ and $(g(X)) \subset A_2$ become equal in the completion $\hat R_2$. Since $A_2$ is a Zariski
ring, this implies that $f(X)$ and $g(X)$ generate the same ideal in $A_2$. Hence $I=(g(X))$ and
because $g(X)$ is homogenous, we have $g(X)= \prod_{j=1} (\gamma_j ' x_1 - \gamma_j x_2)$ for certain $\gamma_j, \gamma_j'\in \Lambda$. Since $I\neq 0$ is prime, 
the polynomial  $ g(X)= \gamma' \cdot x_1 - \gamma\cdot x_2 \neq  0 $ must be linear.
Without restriction of generality we can assume $\gamma\neq 0$. Then $\lambda_1=1, \lambda_2= \gamma'/\gamma$ define the required substitution $\pi$. 

\medskip
For $n\geq 3$ we use induction. If $I\subset (x_2,..,x_n) \subset A_n$,
the assertion is clear. Indeed $m'=(x_2,..,x_n)$ is a prime ideal of $A_{n-1}$
and $A_n/m'$ and $A_n/I$ both have dimension 1, hence $I=m'$. Put $\lambda_1=1$ and $\lambda_i=0$ for $i\geq 2$.
So we can assume that $I$ contains an element $p(X)$ not in $m'$. Up to a scalar in $\Lambda$, then
$p(x_1,0,...,0) = x_1^a + \sum_{i>a} b_i\cdot x_1^i $ for some $a\in \N$ and $b_i\in \Lambda$.
 By the Weierstra\ss\ preparation theorem
we have
$ p(X) = (x_1^a + \sum_{i=0}^{a-1} c_i(x_2,..,x_n)\cdot  x_1^i) \cdot u(X) $
for a unit $u(X)\in A_n$ with coefficients $c_i(X')\in A_{n-1}$ by lemma \ref{WS}. We can therefore assume 
$$ p(X) := x_1^a + \sum_{i=0}^{a-1} c_i(x_2,..,x_n)\cdot  x_1^i  \ \in \ I \ . $$
Intersecting with the subring $A_{n-1} \subset A_n$
analytically generated by $x_2,..,x_n$, we now define the prime ideal $J := I \cap A_{n-1}$ 
in $m' \subset A_{n-1}$. It is a $\varphi$-ideal, for the restriction of $\varphi$ to $A_{n-1}$.
Notice, $A_{n-1}/J \hookrightarrow A_n/I$ is a finite injective ring extension.  Indeed, $A_n$ is finite over the subring $S=\Lambda\{\!\{ p(X), x_2,...,x_n\}\!\}$ by lemma \ref{FINITE} and $S/(S\cap I) \hookrightarrow A_n/I$ is a finite injective  extension of rings. Furthermore $S\cap I =  \Lambda\{\!\{ p(X), x_2,...,x_n\}\!\} \cap I =
p(X)\cdot S + (A_{n-1}\cap I)$, hence $S/(S\cap I) \cong A_{n-1}/J$.
Since for finite injective ring extensions the Krull dimension remains the same,  this implies $\dim(A_{n-1}/J)=1$. 
So, by the induction assumption, the prime ideal $J$ is the kernel of a surjective ring homomorphism $\pi': A_{n-1} \twoheadrightarrow A_1=\Lambda\{\!\{ x\}\!\}$ given by $x_i\mapsto \lambda'_i\cdot x$ for $i=2,..,n$ for some $\lambda'_i$ satisfying the compatibility conditions for $\alpha_2,...,\alpha_n$.
Consider the diagram $$ \xymatrix{    A_n \ar@{->>}[r]^-\pi &  A_2 = \Lambda\{\!\{ x_1,x \}\!\} & \cr
S=\Lambda\{\!\{ p(X), x_2,..,x_n\}\!\}  \ar@{->>}[r] \ar@{^{(}->}[u] &  \Lambda\{\!\{q(x_1,x),x \}\!\}  \ar@{^{(}->}[u]  &   \cr
A_{n-1}\ar@{^{(}->}[u] \ar@{->>}[r]^-{\pi'} &  A_1 = \Lambda\{\!\{ x \}\!\} \ar@{^{(}->}[u] & }  $$
for the corresponding evaluation map $\pi: x_i \mapsto \lambda'_i \cdot x$
for $i=2,..,n$ and $\pi(x_1)=x_1$. Notice, $g(x_1,x) \in \pi(I)$ for $g\in A_2$ holds   iff $g(x_1,x)=f(x_1,\lambda'_2\cdot x,...,\lambda'_n\cdot x)$ holds for some $f(X)\in I$. Hence $g(\alpha_1\cdot x_1, \alpha \cdot x)= f(\alpha_1 \cdot x_1, \alpha\lambda'_2 \cdot x,..., \alpha\lambda'_n \cdot x)
= f(\alpha_1 \cdot x_1, \alpha_2\lambda'_2 \cdot x,..., \alpha_n\lambda'_n \cdot x) = f^\varphi(x_1, \lambda'_2\cdot x,..., \lambda'_n \cdot x)$
is also in $\pi(I)$, using that $f^\varphi \in I$ holds for $f\in I$ and a $\varphi$-ideal $I$. In other words, $\pi(I)$ is a $\varphi$-ideal in $A_2$ where $\varphi$ acts by the eigenvalues $\alpha_1$ resp. $\alpha$ on $x_1$ resp. $x$.
Since $p(X)\in I$ maps to $q(x_1,x)\neq 0$, the image ideal is not zero $\pi(I) \neq 0$. 
A power series $f(X)=\sum_{i=0}^\infty c_i(X')\cdot x_1^i \in A_n$ is in $K:=Kern(\pi: A_n \to A_2)$ if and only if $c_i\in J$
holds for all $i$, hence $K = \widehat{A_n \cdot J}\cap A_n = A_n \cdot J\subseteq I$. 
 We obtain
$$ \xymatrix{ \pi(I)  \ar@{^{(}->}[r]                          &   A_2                           \cr
\pi(I\cap S)=(q)       \ar@{^{(}->}[u] \ar@{^{(}->}[r] &   \pi(S)  \ar@{^{(}->}[u] } $$
and $\pi(I) \cap \pi(S)=\pi((I+K) \cap (S+K))= \pi(I \cap (S+K)) = \pi((I\cap S) + K)= \pi(I\cap S)$.
Since $A_n$ is finite over $S$, $\pi'(A_{n-1}) \cong \pi(S/I\cap S) = \pi(S)/\pi(I\cap S) = \pi(S)/\pi(S)\cap \pi(I) \hookrightarrow A_2/\pi(I) $ is finite and injective again. Since the dimension of $\pi'(A_{n-1})\cong A_1$ is one, the dimension of $A_2/\pi(I)$
is one and  therefore the induction assumption applies. Hence $\pi(I)$ is the kernel of a homomorphism  $\tilde\pi: A_{2} \twoheadrightarrow \Lambda\{\!\{ y \}\!\}$
defined by $x_1 \mapsto \lambda_1\cdot y$ and $x \mapsto \tilde\lambda \cdot y$, satisfying the compatibility condition for $\alpha_1$ and $\alpha$. 
Because the $x_1$-regular power series $q$ is contained 
in $\pi(I)$, this implies $\tilde\lambda \neq 0$.
The composition
$$ \xymatrix{ A_{n} \ar@{->>}[r]^\pi &  A_{2}  \ar@{->>}[r]^-{\pi'} & \Lambda\{\!\{ y \}\!\} \cr } $$
defines the desired specialization $x_1 \mapsto \lambda_1\cdot y$ and   
$x_i \mapsto \lambda_i  \cdot y$ with $\lambda_i = \tilde\lambda\lambda'_i$ for $i=2,...,n$. It is easy to see that the compatibility conditions for the $\alpha_1,..,\alpha_n$ are satisfied. For $\lambda_1=0$ there is nothing to show,
and $\lambda_1\neq 0$ implies $\alpha_1=\alpha$. Since $\alpha_i=\alpha$ for all $i\geq 2$ with
$\lambda_i\neq 0$, this completes the proof of the one-dimensional case of proposition \ref{P}. 

\medskip
For the higher dimensional situation, we use both induction on $n$ and then for fixed $n$ induction on $r=\dim(A_n/I)$.
Suppose $I\neq m$ is a prime $\varphi$-ideal in $A_n$, now of dimension $r=\dim(A_n/I) >1$. 
Consider the intersection of $Spec(A_n/I)$ with the hypersurface 
$Spec(A_n/P)\cong A_{n-1}$ defined by the $\varphi$-ideal $P=(x_1)\subset A_n$.
Then either  $P \subseteq I$, i.e.  $Spec(A_n/I)$ is contained in the hypersurface $Spec(A_n/P)$, and we may replace 
$A_n$ by $A_n/P=A_{n-1}$ and $I$ by $I/P$ to conclude by induction on $n$ using $\dim ((A_n/P)/(I/P))=\dim(A_n/I)\neq 0$.
Or the intersection of $Spec(A_n/I)$ with $Spec(A_n/P)$, given by the $\varphi$-ideal $(I,x_1) \subseteq m \subset  A_n$,
has the property $I \neq (I,x_1)$. Since $I$ is prime, $A_n/I$ is a domain. The image $\overline x_1$ of $x_1$ in $A_n/I$ therefore  is not a zero divisor nor a unit. Hence $\overline x_1$ is not a zero divisor of $A_n/I$ nor a unit.
Since $A_n/I$ is a local and noetherian, therefore $\dim(A_n/(I,x_1))=\dim((A_n/I)/ \overline x_1 (A_n/I)) = \dim(A_n/I) - 1$ by [AM, 11.18].
We then find a minimal prime ideal of $A_n/(I,x_1)$ whose
preimage $I'$ in $A_n$ 
satisfies $dim(A_n/I') = \dim(A_n/I)-1 = r-1 >0$, hence $I' \neq m$. Since the prime ideal $I'$ defines one of the finitely many irreducible components of the spectrum of the $\varphi$-ideal $(I,x_1)$, replacing $\varphi$ with some power we can assume that $I'$ is a $\varphi$-ideal of $A_n$.  
Hence there exists $\pi': A_n/I' \twoheadrightarrow A_1$ by the induction assumption on $r$, so that finally $\pi: A_n/I \twoheadrightarrow A_n/I' \twoheadrightarrow A_1$ defines the desired
evaluation morphism. \qed

\bigskip\noindent

\medskip
\goodbreak

\centerline {\bf Appendix}

\bigskip\noindent

\medskip  
Let $A_n=\Lambda\{\!\{t_1,...,t_n\}\!\}$ denote the ring of power series in the variables $t_1,...,t_n$ with coefficients in some finite extension field $E_\lambda$ of $\Q_\ell$ that are convergent in $\Lambda$ with some positive radius $r$ of convergence with respect to the nonarchimedean norm on $\Lambda$, with $r$ and $E_\lambda$ depending on the power series. 

\begin{Lemma} \label{WS} The Weierstra\ss\ preparation theorem holds
for $A_n$. The ring $A_n$ is a regular noetherian local ring of Krull dimension $n$,
hence $A_n$ is a normal factorial domain. \end{Lemma}

\medskip
{\it Proof}. Any substitution $t_i \mapsto \lambda_i \cdot t_i$ for $\lambda_i\in \Lambda^*$ defines an automorphism of $A_n$. For the proof of the Weierstra\ss\ preparation theorem  suppose given $G\in A_n$ and a $t_1$-regular $F\in A_n$ such that $F(t_1,0,..,0) =c\cdot t_1^a$ plus terms of higher order. 
We have to show the existence of $U\in A_n$ and $R_0,..,R_{a-1}\in A_{n-1}=\Lambda\{\!\{ t_2,...,t_n\}\!\}=:\Lambda\{\!\{ T'\}\!\}$ such that
$G(T) = U(T)\cdot F(T) + \sum_{i=0}^{a-1} R_i(T') \cdot t_1^i$. We may assume that $F$ is not a unit in $A_n$ and that
$G(0,..,0)$ has absolute value $\leq 1$. Then, by a suitable substitution $T \mapsto \lambda \cdot T$, we can assume $F,G \in {\frak o}_\lambda[[t_1,...,t_2]] \subset A_n$ for some subring ${\frak o}_\lambda$ of $\Lambda$ that is finite over $\Z_\ell$ with maximal ideal $m_\lambda$. Replacing $F$ by $c^{-1-a}F(ct_1,c^{a+1}t_2,..,c^{a+1}t_n)$ for some nonzero $c\in {\frak o}_\lambda$, we may assume  $F(t_1,0,..,0)=t_1^a + \sum_{i>a} c_i(T')\cdot t_1^i$
for $c_i(T') \in A_{n-1}$ and $c^{-1-a}F(ct_1,c^{a+1}t_2,..,c^{a+1}t_n)\in {\frak o}_\lambda[[t_1,...,t_n]] $. 
Again, by replacing $F(T)$ with $b^{-a}F(bt_1,b^at_2,...,b^at_n)$ for some $b\in m_\lambda$, we can assume $F(T) \equiv t_1^a $ modulo $(m_\lambda,t_2,...,t_n)$. 
Since $G(bct_1,c(bc)^at_2,..,c(bc)^at_n)$
is in ${\frak o}_\lambda[[t_1,...,t_n]]$ and $b^{a}c^{1+a}$ is a unit in $A_n$, the proof of the preparation theorem is reduced to  [GL, prop. A.2.1(i)], i.e. the assertion that the Weierstra\ss\ preparation theorem holds
for  ${\frak o}_\lambda[[t_1,...,t_n]]$.  
Since, up to a linear coordinate change, any nontrivial ideal $I$ of $A_n$ contains a $t_1$-regular element
$F(T)=t_1^a + c_1(T')t_1^{a-1} + ... + c_a(T')$ with coefficients $c_\nu \in \Lambda\{\!\{t_2,...,t_n\}\!\} \cong A_{n-1}$,
by the Weierstra\ss\ preparation theorem then $F$ and $I'= I \cap \{ A_{n-1} t_1^{a-1} + ... + A_{n-1}\}$ generate
$I$. To show that $A_n$ is noetherian we can assume $A_{n-1}$ to be noetherian by induction, so $I'$ is a finitely generated $A_{n-1}$-module and its generators together with $F$ generate $I$ as an $A_n$-module. This proves that $A_n$ is noetherian.
It is easy to see that for any powers series $F(T)\in A_n$ with $F(0) \neq 0$ the formal power series
$1/F(T)$ again has positive radius of convergency. Hence $A_n$ is a local ring with maximal ideal $m=(t_1,...,t_n)$.
Since $\hat A_n$ is isomorphic to the regular  ring $\hat R_n$
of formal powers series over $\Lambda$, $\hat A_n$ is a regular local ring. \qed

\medskip
The regular noetherian ring
$R_n = \Lambda \otimes_{\Z_\ell} {\Z}_\ell[[t_1,....,t_{n}]]$
is a subring  of the ring $A_n$. The 
completions $\hat A_n$ resp. $\hat R_n$ of $A_n$ resp. $R_n$ with respect to $m=(t_1,..,t_n)$ coincide with the formal power series ring $\Lambda[[t_1,...,t_n]]$. Any power series $P(t_1,..,t_n)\in R_n$ with $P(0,...,,0)\neq 0$
is a unit in $A_n$. Hence the localisation $\tilde R_n=(R_n)_{m} $ of $R_n$ with respect to its maximal ideal $m=(t_1,..,t_n)$ is a local noetherian subring of $A_n$. Both local rings $\tilde R_n$ and $A_n$ have the same completion (for their maximal ideals), namely the formal power series ring $\hat R_n$. 
The ideals $I$ in $m \subset R_n$ correspond one-to-one to ideals in $\tilde R_n$. Since Noetherian local rings are Zariski rings [ZS, p.264], any ideal $I$ of a noetherian local
ring $R$  can be recovered from its completion $\hat I$ by intersection $I= \hat I \cap R$, and furthermore $\hat I = I \cdot \hat R$ holds. See  [ZS, VIII, $\S 2$, thm. 5, cor. 2] resp. [ZS, VIII, $\S 4$, thm. 8] and for the definition of Zariski rings [ZS, VIII $\S$ 4]. For $I\subset R_n$ this applies for both $\tilde R_n \cdot I\subset \tilde R_n$ and $A_n\cdot I \subset A_n$, hence $I$ can be recovered from $\tilde R_n \cdot I$, which can be recovered from $\hat R_n \cdot I$ or $A_n\cdot I$. 
In particular, for any ideal $I\subset R_n$ contained in $(t_1,...,t_n)$ the ideal
$\tilde I = I \cdot A_n = \hat I \cap A_n$ generated by $I$ in $A_n$  is maximal resp. zero if and only if $\tilde I$ is maximal resp. zero
in $A_n$. Similarly, two ideals $J,J'$ in $A_n$ are equal iff $J\cdot \hat R_n$ and $J'\cdot \hat R_n$ are equal.

\begin{Lemma} \label{comm}
A normal noetherian domain $R$ is factorial if and only if every prime ideal $I$ of height $ht(I)=1$  is a
principal ideal. 
\end{Lemma}

{\it Proof}. [Bourbaki 7.3, no. 2, thm. 1]
or [Sa, thm. 3.2,5.3, p.6 cor.]. 

\begin{Lemma} \label{FINITE}
Let $p(X)$ be a $x_1$-regular homogenous polynomial $x_1^a + \sum_{\nu< a} c_\nu(X') \cdot x_1^\nu$ of degree $a>0$ in $A_n=\Lambda\{\!\{ x_1,x_2,...,x_n\}\!\}$ with coefficients $c_\nu(X')$ in $A_{n-1}=\Lambda\{\!\{x_2,...,x_n\}\!\}$ for $i=0,...,a-1$ and with $c_0(X')$ in $(x_2,...,x_n)$.
Then $A_n$ is a finite ring extension of its subring $\Lambda\{\!\{ p(X),x_2,...,x_n\}\!\}$.
\end{Lemma}

{\it Proof}. Any $g(X)$ in $\Lambda\{\!\{ X\}\!\} $ can be written in the form $g(X)=u(X) \cdot p(X) +  \sum_{i=0}^{a-1} r_i(X') \cdot x_1^i$ by the Weierstra\ss\ preparation theorem (lemma \ref{WS}). If we apply this iteratively for $u(X)$ instead of $g(X)$ and continue, we obtain formal power series $f_i(y_1,....,y_n)\in \Lambda[[y_1,...,y_n]]$ so that
$g(X) = \sum_{i=0}^{a-1} f_i(p(X),x_2,...,x_n)\cdot x_1^i$ holds in $\Lambda[[x_1,...,x_n]]$. 
To prove
$\Lambda\{\!\{ X\}\!\} = x_1^{a-1} \cdot \Lambda\{\!\{ p,X'\}\!\} +  x_1^{a-2}\Lambda\{\!\{ p,X'\}\!\}  + \cdots  +  \Lambda\{\!\{ p,X'\}\!\}$ it suffices to show $f_i\in \Lambda\{\!\{ Y \}\!\}$. If $g(X)$ and $p(X)$ both have integral coefficients in ${\frak o}_\lambda$ 
with $p(X)\equiv x_1^a $ modulo $(m_\lambda,x_2,...,x_n)$,  then all $f_i(Y)$ have integral coefficients by the Weierstra\ss\ preparation theorem for the ring ${\frak o}_\lambda[[x_1,...,x_n]]$. The general case can be easily reduced to this by the method used in the proof of lemma \ref{WS}. We may assume $g(0)=0$ and replace $g(X)$ by $g(c \cdot x_1,c^a \cdot x_2,...,c^a \cdot x_n)$ and $p(X)$ by $\tilde p(X)=c^{-a} \cdot p(c \cdot x_1, c^a \cdot x_2,..., c^a x_n)$ to show  $f_i(y_1,...,y_n)= c^{-i} \cdot g_i(c^{-a}\cdot y_1, y_2,..,y_n)$ for certain $g_i\in {\frak o}_\lambda[[x_1,...,x_n]]$. Hence all $f_i(y_1,..,y_n)$ are locally convergent. 
\qed

\bigskip\noindent
\bigskip\noindent

\goodbreak

\bigskip\noindent
\centerline{\bf References}.

\bigskip\noindent

\medskip\noindent
[AS] Abbes A., Saito T., {\it Ramification and cleanliness}, 
Tohoku Mathematical Journal, Centennial Issue, 63 No. 4 (2011), 775-853.

\medskip\noindent
[AM] Atiyah M.F., Macdonald I.G., {\it Introduction to commutative algebra}, Addison Wesley (1969)

\medskip\noindent
[BBD] Beilinson A. A., Bernstein J., Deligne P., {\it Faiscaux pervers}, in: analyse et topologie sur les espaces singuliers (I), asterisque 100, SMF (1982) 

\medskip\noindent
[E] Ekedahl T., {\it On the adic formalism}, in: {\it The Grothendieck Festschrift Vol. II}, Progr.  Math. 87, Birkh\"auer, Boston, (1990), 197 - 218

\bigskip\noindent
[D] Deligne P., {\it La conjecture de Weil, II.}, Inst. Hautes Etudes Sci. Publ. Math. 52 (1980), 137 - 252

\bigskip\noindent
[D2] Deligne P., {\it Finitude de l'extensions de $\mathbb Q$ engendree par les traces de Frobenius, en characterisque finie}, Moscow Mathematical Journal 12 (2012)

\medskip\noindent
[Dr] Drinfeld V., {\it On a conjecture of Deligne}, Moscow Math. Journal 12 (2012)

\medskip\noindent
[Dr2] Drinfeld V., {\it On a conjecture of Kashiwara}, Math. Res. Lett. 8 (2001) 713 - 728

\medskip\noindent
[FK] Franecki J., Kapranov M., {\it The Gauss map and a noncompact Riemann-Roch formula for constructible sheaves 
on semiabelian varieties}, Duke Math. J. 104 no.1 (2000), 171 - 180

\medskip\noindent
[GL] Gabber O., Loeser F., {\it Faiscaux pervers $\ell$-adiques sur un tore}, Duke Math. Journal vol. 83, no. 3 (1996),
501 - 606


\medskip\noindent
[KrW] Kr\"amer T., Weissauer R., {\it Vanishing Theorems for constructible sheaves on abelian varieties}, arXiv: 1111.6095v5; to appear in 
J. Alg. Geom.

\medskip\noindent
[KrW2] Kr\"amer T., Weissauer R., {\it On the Tannaka group attached to the Theta divisor of a generic principally polarized abelian variety}, arXiv: 1309.3754; to appear in 
Math. Zeitschrift.

\medskip\noindent
[L] Laumon G., {\it Letter to Gabber and Loeser} (22/12/91)

\medskip\noindent
[Sa] Samuel P., {\it Lectures
On Unique Factorization Domains}, Tata Institute, Bombay (1964)

\medskip\noindent 
[S] Serre J.P., {\it Algebraic groups and class fields}, Springer (1988)

\medskip\noindent
[W] Weissauer R., {\it On the rigidity of BN-sheaves}, arXiv: 1204.1929 

\medskip\noindent
[W2] Weissauer R., {\it Why certain Tannaka groups attached to abelian varieties are almost connected}, arXiv: 1207.4039 

\medskip\noindent
[W3] Weissauer R., {\it Degenerate Perverse Sheaves on Abelian Varieties}, arXiv: 1204.2247

\medskip\noindent
[ZS] Zariski O., Samuel P., {\it Commutative Algebra, Volume II}, Springer (1960)
\end{document}